\documentclass{amsart}

\usepackage{amsfonts}
\usepackage{amsmath}
\usepackage{amsthm}
\usepackage{amssymb}
\usepackage[english]{babel}
\usepackage{graphics}
\usepackage{latexsym}
\usepackage{longtable} \usepackage{mathrsfs}
\usepackage{graphicx}
\usepackage{wrapfig} 
\usepackage[pdftex,colorlinks=true,linkcolor=blue,citecolor=cyan]{hyperref}

\newtheorem{theorem}{Theorem}
\newtheorem{corollary}[theorem]{Corollary}
\newtheorem{lemma}[theorem]{Lemma}

\newtheorem{claim}[theorem]{Claim}
\newtheorem{example}[theorem]{Example}

\theoremstyle{definition}
\newtheorem{definition}[theorem]{Definition}

\newcommand{\msL}{\mathscr{L}}

\newcommand{\mD}{\mathcal{D}}

\newcommand{\R}{\mathbb{R}}
\newcommand{\N}{\mathbb{N}}

\newcommand{\D}{\mathrm{D}}

\newcommand{\ms}{\medskip}

\newcommand{\de}{\delta}

\newcommand{\e}{\varepsilon}

\newcommand{\Om}{\Omega}
\newcommand{\om}{\omega}

\newcommand{\weak }{\, -\!\!\!\!-\!\!\!\rightharpoonup}

\newcommand{\larrow}{\longrightarrow}
\newcommand{\ot}{\otimes}

\newcommand{\lmapsto}{\longmapsto}
\newcommand{\ri}{\rightarrow}

\newcommand{\sub}{\subseteq}

\newcommand{\by}{\times}

\newcommand{\sgn}{\mathrm{sgn}}
\newcommand{\ess}{\mathrm{ess}}

\newcommand{\bt}{\begin{theorem}}\newcommand{\et}{\end{theorem}}
\newcommand{\bd}{\begin{definition}}\newcommand{\ed}{\end{definition}}
\newcommand{\bl}{\begin{lemma}}\newcommand{\el}{\end{lemma}}
\newcommand{\beq}{\begin{equation}}\newcommand{\eeq}{\end{equation}}
\newcommand{\bc}{\begin{claim}}\newcommand{\ec}{\end{claim}}
\newcommand{\bex}{\begin{example}}\newcommand{\eex}{\end{example}}
\newcommand{\bcor}{\begin{corollary}}\newcommand{\ecor}{\end{corollary}}
\newcommand{\bp}{\begin{proof}}\newcommand{\ep}{\end{proof}}

\newcommand{\BPT}{\medskip \noindent \textbf{Proof of Theorem} }

\numberwithin{equation}{section}

\begin{document}

\title[Vectorial Absolute Minimisers under Minimal Assumptions]{Existence of $1D$ Vectorial Absolute Minimisers in $L^\infty$ under Minimal Assumptions}

\author{Hussien Abugirda}
\address{Department of Mathematics , College of Science, University of Basra, Basra, Iraq $\&$ Department of Mathematics and Statistics, University of Reading, Whiteknights, PO Box 220, Reading RG6 6AX, UK}
\email{h.a.h.abugirda@student.reading.ac.uk}

\author{Nikos Katzourakis}

\address{Department of Mathematics and Statistics, University of Reading, Whiteknights, PO Box 220, Reading RG6 6AX, UK}
\email{ n.katzourakis@reading.ac.uk}


\date{}



\begin{abstract} We prove the existence of vectorial Absolute Minimisers in the sense of Aronsson to the supremal functional $E_\infty(u,\Omega') = \|\mathscr{L}(\cdot,u,\D u)\|_{L^\infty(\Omega')}$, $\Omega'\Subset \Omega$, applied to $W^{1,\infty}$ maps $u:\Om\subseteq \R\longrightarrow \R^N$ with given boundary values. The assumptions on $\msL$ are minimal, improving earlier existence results previously established by Barron-Jensen-Wang and by the second author.  
\end{abstract}

\maketitle


\section{Introduction} \label{section1}

The main goal of this paper is to prove the existence of vectorial Absolute Minimisers with given boundary values to the supremal functional 
\beq \label{1.1}
E_\infty (u,\Om')\, :=\, \underset{x\in \Om'}{\ess\,\sup} \, \mathscr{L}\left(x,u(x),\D u(x) \right), \quad  \, u\in W^{1,\infty}_{\text{loc}}(\Om,\R^N),\ \Om'\Subset \Om,
\eeq
applied to maps $u: \Om\sub \R \larrow \R^N$, $N\in \N$, where $\Om$ is an open interval and $\msL \in C(\Om \by \R^N \by  \R^N )$ is a non-negative continuous function which we call Lagrangian and whose arguments will be denoted by $(x,\eta,P)$. By Absolute Minimiser we mean a map $u\in W^{1,\infty}_{\text{loc}}(\Om,\R^N)$ such that
\beq \label{1.2}
\ \ \ E_\infty (u,\Om')\,\leq\, E_\infty (u+\phi,\Om'),
\eeq
for all $\Om'\Subset \Om$ and all $\phi\in W^{1,\infty}_0(\Om',\R^N)$. This is the appropriate minimality notion for supremal functionals of the form \eqref{1.1}; requiring at the outset minimality on all subdomains is necessary because of the lack of additivity in the domain argument. The study of \eqref{1.1} was pioneered by Aronsson in the 1960s \cite{A1}-\cite{A5} who considered the case $N=1$. Since then, the (higher dimensional) scalar case of $u :\Om\sub \R^n \ri\R$ has developed massively and there is a vast literature on the topic (see for instance the lecture notes \cite{C, K7}). In the case the Lagrangian is $C^1$, of particular interest has been the study of the (single) equation associated to \eqref{1.1}, which is the equivalent of the Euler-Lagrange equation for supremal functionals and is known as the ``Aronsson equation":
\beq \label{Ar}
A_\infty u \, :=\, \D\big(\msL(\cdot,u,\D u)\big)\msL_P(\cdot,u,\D u)\,=\,0.
\eeq
In \eqref{Ar} above, the subscript denotes the gradient of $\msL(x,\eta,P)$ with respect to $P$ and, as it is customary, the equation is written for smooth solutions. Herein we are interested in the vectorial case $N\geq 2$ but in one spatial dimension. Unlike the scalar case, the literature for $N\geq 2$ is much more sparse and starts much more recently. Perhaps the first most important contributions were by Barron-Jensen-Wang \cite{BJW1, BJW2} who among other deep results proved the existence of Absolute Minimisers for \eqref{1.1} under certain assumptions on $\msL$ which we recall later. However, their contributions were at the level of the functional and the appropriate (non-obvious) vectorial analogue of the Aronsson equation was not known at the time. The systematic study of the vectorial case of \eqref{1.1} (actually in the general case of maps $u:\Om\sub\R^n \larrow \R^N$) together with its associated system of equations begun in the early 2010s by the second author in a series of papers, see \cite{K1}-\cite{K6}, \cite{K8}-\cite{K12} (and also the joint contributions with Croce, Pisante and Pryer \cite{CKP, KP1, KP2}). The ODE system associated to \eqref{1.1} for smooth maps $u:\Om\sub\R \larrow \R^N$ turns out to be
\beq \label{1.3a}
\mathcal{F}_\infty \big(\cdot,u,\D u,\D^2u\big)\,=\,0, \quad \text{on }\Om,
\eeq
where
\beq  \label{1.4a}
\begin{split}
\mathcal{F}_\infty (x,\eta,P,X) \,  :=&  \, \Big[  \msL_P(x,\eta,P) \ot \msL_P(x,\eta,P)\\
 &  +  \, \msL(x,\eta,P) [\msL_P(x,\eta,P)]^\bot  \msL_{PP}(x,\eta,P)\Big]X  
\\ 
 &+\,  \Big( \msL_\eta(x,\eta,P) \cdot P \, +\, \msL_x(x,\eta,P) \Big) \msL_P  (x,\eta,P)\\
  & + \,  \msL  (x,\eta,P)\big[\msL_P(x,\eta,P)\big]^\bot \Big(\msL_{P\eta}(x,\eta,P)P
 \\
 & +  \,  \msL_{Px}  (x,\eta,P) \,-\, \msL_\eta(x,\eta,P) \Big).
\end{split}
\eeq
Quite unexpectedly, in the case $N\geq 2$ the Lagrangian needs to be $C^2$ for the equation to make sense, whilst the coefficients of the full system are \emph{discontinuous}; for more details we refer to the papers cited above. In \eqref{1.4a} the notation of subscripts symbolises derivatives with respect to the respective variables and $\smash{\big[\msL_P(x,\eta,P)\big]}^\bot$ is the orthogonal projection to the hyperplane normal to $\msL_P(x,\eta,P) \in \R^N$:
\beq  \label{1.5a}
\ \ \big[\msL_P(x,\eta,P)\big]^\bot\, :=\, \mathrm{I} - \, \sgn\big(\msL_P(x,\eta,P)\big) \ot \sgn\big(\msL_P(x,\eta,P)\big).
\eeq
The system \eqref{1.3a} reduces to the equation \eqref{Ar} when $N=1$. In the paper \cite{K9} the existence of an absolutely minimising \emph{generalised solution} to \eqref{1.3a} was proved, together with extra partial regularity and approximation properties. Since  \eqref{1.3a} is a quasilinear non-divergence degenerate system with discontinuous coefficients, a notion of appropriately defined ``weak solution" is necessary because in general solutions are non-smooth. To this end, the general new approach of $\mD$-solutions which has recently been proposed in \cite{K8} has proven to be the appropriate setting for vectorial Calculus of Variations in $L^\infty$ (see \cite{K8}-\cite{K10}), replacing to some extent viscosity solutions which essentially apply only in the scalar case.

Herein we are concerned with the existence of absolute minimisers to \eqref{1.1} without drawing any connections to the differential system \eqref{1.3a}. Instead, we are interested in obtaining existence under the weakest possible assumptions. Accordingly, we establish the following result.

\begin{theorem}  \label{theorem1}
Let $\Om\sub \R$ be a bounded open interval and let also
\[
\begin{split}
\mathscr{L}\ :\ \overline{\Om}\by\R^N \by\R^N\larrow [0,\infty),\ \ \ \
\end{split}
\]
be given continuous function with $N\in \N$. We assume that:
\begin{enumerate}
\item For each $ (x,\eta)  \in  \overline{\Om}\by\R^N$, the function $P \lmapsto \mathscr{L}(x,\eta,P) $ is level-convex, that is for each $t\geq 0$ the sublevel set
\[
\big\{P\in \R^N \,:\, \mathscr{L}(x,\eta,P) \leq t\big\}
\]
is a convex set in $\R^N$.

\item there exist non-negative constants  $ C_{1}, C_{2}, C_{3}$, and $ 0<q\leq r<+\infty$ and a positive locally bounded function $h:\R\by\R^N\larrow [0,+\infty)$ such that for all $(x,\eta,P) \in \overline{\Om}\by\R^N \by\R^N$ 
\[
\ \ \ \ C_{1}|P|^q \,-\, C_{2} \, \leq \, \mathscr{L}(x,\eta,P)\, \leq  \, h(x,\eta)|P|^r \,+\, C_{3}.
\]
\end{enumerate}
Then, for any affine map  $b : \R\larrow \R^N$, there exist a  vectorial Absolute Minimiser $u^\infty\in W^{1,\infty}_b(\Om,\R^N)$ of the supremal functional \eqref{1.1} (Definition \eqref{1.2}).
\end{theorem}

Theorem \ref{theorem1} generalises two respective results in the both the papers \cite{BJW1} and \cite{K9}. On the one hand, \emph{in \cite{BJW1} Theorem \ref{theorem1} was established under the extra assumption $C_2=C_3=0$ which forces $\msL(x,\eta,0)=0$, for all $(x,\eta)\in\R \by \R^N$.} Unfortunately this requirement is incompatible with important applications of \eqref{1.1} to problems of $L^\infty$-modelling of variational Data Assimilation (\emph{4DVar}) arising in the Earth Sciences and especially in Meteorology (see \cite{B, BS, K9}). An explicit model of $\msL$ is given by
\beq \label{1.4}
\mathscr{L}(x,\eta,P)\, :=\, \big|k(x)-K(\eta) \big|^2 +\,   \big|P-\mathscr{V}(x,\eta) \big|^2,
\eeq
and describes the ``error" in the following sense: consider the problem of finding the  solution $u$ to the following ODE coupled by a pointwise constraint:
\[
\ \ \ \D u(t) = \mathscr{V}\big(t,u(t)\big) \ \ \ \& \ \ \ K(u(t)) = k(t), \ \ \ t\in \Om.
\]
Here  $\mathscr{V} : \Om \by \R^N\larrow \R^N$ is a time-dependent vector field describing the law of motion of a body moving along the orbit described by  $u:\Om\sub \R \larrow \R^N$ (e.g.\ Newtonian forces, Galerkin approximation of the Euler equations, etc), $k: \Om\sub \R \larrow \R^M$ is some partial ``measurements" in continuous time along the orbit and $K: \R^N \larrow \R^M$ is a submersion which corresponds to some component of the orbit that is observed. We interpret the problem as that $u$ should satisfy the law of motion and also be compatible with the measurements along the orbit. Then minimisation of \eqref{1.1} with $\msL$ as given by \eqref{1.4} leads to a uniformly optimal approximate solution without ``spikes" of large deviation of the prediction from the actual orbit.

On the other hand, in the paper \cite{K9} Theorem \ref{theorem1} was proved under assumptions allowing to model Data Assimilation but strong convexity, smoothness and structural assumptions  were imposed, allowing to obtain stronger results accordingly. In particular, the Lagrangian was assumed to be radial in $P$, which means it can be written in the form
\[
\ \ \ \ \ \msL(x,\eta,P)\, :=\,\mathscr{H}\Big(x,\eta,\frac{1}{2} \big|P-\mathscr{V}(x,\eta) \big|^2\Big).
\]
In this paper we relax the hypotheses of both the aforementioned results.

\section{Proof of the main result}

In this section we establish Theorem \ref{theorem1}. In the proof we will utilise the following lemma essentially proved in the paper \cite{BJW1} which we recall right below for the convenience of the reader.

\begin{lemma} \label{lemma1}\
\big(cf. Lemma 2.2 in \cite{BJW1}\big)
In the setting of theorem \ref{theorem1} and under the same hypotheses, for a fixed affine map $b : \R\larrow \R^N$,  set 
\[
\begin{split}
C_{m}\, &:=\, \inf\Big\{E_{m}(u,\Om) \,:\, u\in W^{1,qm}_b(\Om,\R^N)\Big\},
\\
C_\infty \, &:=\, \inf\Big\{E_\infty(u,\Om)\, :\, u\in W^{1,\infty}_b(\Om,\R^N)\Big\}.
\end{split}
\]
where $E_\infty$ is as in  \eqref{1.1} and
\beq \label{2.1a}
E_m(u,\Om)\,:=\,\int_\Om \msL\big(x,u(x),\D u(x)\big)^m\,dx.
\eeq
Then, there exist  $u^\infty\in W^{1,\infty}_b(\Om,\R^N)$ which is a (mere) minimiser of \eqref{1.1} over $W^{1,\infty}_b(\Om,\R^N)$ and a sequence of approximate minimisers $\{u^m\}^\infty_{m=1}$ of \eqref{2.1a} in the spaces $W^{1,qm}_b(\Om,\R^N)$ such that, for any $s\geq 1$, 
\[
\text{$u^m \weak u^\infty$, weakly as $m\ri \infty$ in $W^{1,s}(\Om,\R^N)$} 
\]
along a subsequence. Moreover,
\beq  \label{2.1}
E_\infty(u^\infty,\Om)\, =\,  C_\infty\, =\, \underset{m\ri\infty}{\lim}(C_m)^{\frac{1}{m}} .
\eeq
By approximate minimiser we mean that $u^m$ satisfies
\beq  \label{2.2}
\big| E_m(u^m,\Om) - C_m \big| \,<\, 2^{-m^2},
\eeq
Finally, for any $A\sub\Om$ measurable of positive measure the following lower semicontinuity inequality holds
\beq  \label{2.3}
E_\infty(u^\infty,A)\, \leq\, \underset{m\ri\infty}{\liminf}\, E_m(u^m,A)^{\frac{1}{m}} .
\eeq
\end{lemma}

The idea of the proof of \eqref{2.2} is based on the use of Young measures in order to bypass the lack of convexity for the approximating $L^m$ minimisation problems (recall that $\msL(x,\eta,\cdot)$ is only assumed to be level-convex); without weak lower-semicontinuity of $E_m$, the relevant infima of the approximating functionals may not be realised. For details we refer to \cite{BJW1} (this method of \cite{BJW1} has most recently been applied to higher order $L^\infty$ problems, see \cite{KP2}). We also note that \eqref{2.3} has been established in p.\ 264 of \cite{BJW1} in slightly different guises, whilst the scaling of the functionals $E_m$ is also slightly different therein. However, it is completely trivial for the reader to check that their proofs clearly establish our Lemma \ref{lemma1}.

\BPT \ref{theorem1}. Our goal now is to prove that the candidate $u^\infty$ of Lemma \ref{lemma1} above is actually an Absolute Minimiser of \eqref{1.1}, which means we need to prove $u^\infty$ satisfies \eqref{1.2}. 

The method we utilise follows similar lines to those of \cite{K9}, although technically has been slightly simplified. The main difference is that due to the weaker assumptions than those of \cite{K9}, we invoke the general Jensen's inequality for level-convex functions. In \cite{K9} the Lagrangian was assumed to be radial in the third argument, a condition necessary and sufficient for the symmetry of the coefficient matrix multiplying the second derivatives in \eqref{1.3a}; this special structure of $\msL$ led to some technical complications. Also, herein we have reduced the number of auxiliary parameters in the energy comparison map (defined below) by invoking a diagonal argument.

Let us fix  $\Om'\Subset\Om$. Since $\Om'$ is a countable disjoint union of open intervals, there is no loss of generality in assuming $\Om'$ is itself an open interval. By a simple rescaling argument, it suffices to assume that $\Om'=(0,1) \sub \R$. For, let $\phi\in W^{1,\infty}_0((0,1),\R^N)$ be an arbitrary variation and set $\psi^\infty := u^\infty +\phi$. In order to conclude, it suffices to establish
\[
E_\infty \big(u^\infty,(0,1)\big) \, \leq\, E_\infty \big(\psi^\infty,(0,1)\big).
\]
Obviously, $u^\infty(0)=\psi^\infty(0)$ and $u^\infty(1)=\psi^\infty(1)$. We define the energy comparison function $\psi^{m,\de}$, for any fixed $0<\de<1/3$ as 
\[
\psi^{m,\de}(x)\,:=\, 
\left\{
\begin{array}{ll}
\left(\dfrac{\de-x}{\de}\right) u^m(0) + \left(\dfrac{x}{\de}\right)\psi^\infty(\de), &   x\in (0,\de], \ms\\
\ \psi^\infty(x), & x\in (\de,1-\de), \ms\\
\left(\dfrac{1-x}{\de} \right)\psi^\infty(1-\de) + \left(\dfrac{x-(1-\de)}{\de}\right)u^m(1), &   x\in [1-\de,1) ,
\end{array}
\right.
\]
where $m\in \N \cup\{\infty\}$. Then, $\psi^{m,\de}-{u^m} \in W^{1,\infty}_0\big((0,1),\R^N\big)$ and 
\[
\D \psi^{m,\de}(x)\,=\, 
\left\{
\begin{array}{ll}
 \dfrac{\psi^\infty(\de) - u^m(0)}{\de}, &  \text{on } (0,\de), \ms\\
 \D \psi^\infty, &  \text{on } (\de,1-\de), \ms \ms\\
 \dfrac{ \psi^\infty(1-\de) - u^m(1) }{-\de}, &   \text{on } (1-\de,1) .
\end{array}
\right.
\]
Now, note that 
\beq \label{2.4}
\ \ \ \ \text{$\psi^{m,\de} \larrow \psi^{\infty,\de}$ in $W^{1,\infty}\big((0,1),\R^N\big)$, as $m\ri \infty$,}
\eeq
because $\psi^{m,\de} \larrow \psi^{\infty,\de}$ in $L^\infty\big((0,1),\R^N\big)$ and for a.e.\ $x\in (0,1)$ we have
\[
\begin{split}
\Big| \D\psi^{m,\de}(x) - \D\psi^{\infty,\de}(x)\Big|\, &=\, \chi_{(0,\de)}\frac{|u^\infty(0)-u^m(0)|}{\de} \,+\, \chi_{(1-\de,1)}\frac{|u^\infty(1)-u^m(1)|}{\de}
\\
& \leq\, \left(\frac{1}{\de} +\frac{1}{\de}\right) \|u^m-u^\infty\|_{L^\infty(\Om)}
\\
&  =\, o(1),
\end{split}
\]
as $m\ri \infty$ along a subsequence. 

Now, recall that $\psi^{m,\de}=u^m$ at the endpoints $\{0,1\}$. Let us also remind to the reader that after the rescaling simplification, $(0,1)$ is a subinterval of $\Om \sub \R$ whilst \eqref{2.2} holds only for the whole of $\Om$. Since $u^m$ is an approximate minimiser of \eqref{2.1a} over $W^{1,m}_b(\Om,\R^N)$ for each $m\in \N$, by utilising the approximate minimality of $u^m$ (given by \eqref{2.2}), the additivity of $E_m$ with respect to its second argument, we obtain the estimate 
\[
E_m \big(u^m,(0,1)\big)  \, \leq \, E_m \big(\psi^{m,\de} ,(0,1)\big)  \,+\, 2^{-m^2}.
\]
Hence, by H\"older inequality
\beq \label{2.5}
\begin{split}
E_m \big(u^m,(0,1)\big)^{\frac{1}{m}} \, &\leq \, E_m \big(\psi^{m,\de} ,(0,1)\big)^{\frac{1}{m}} \,+\, 2^{-m}
\\
&\leq \, E_\infty \big(\psi^{m,\de},(0,1)\big) \,+\,  2^{-m}.
\end{split} 
\eeq
On the other hand, we have 
\[
\begin{split}
E_\infty \big(\psi^{m,\de},(0,1)\big)\, =& \, \max\Big\{ E_\infty  \big(\psi^{m,\de},
(0,\de)\big), \\
&\ \ \ \ \ \  \ \ \ E_\infty  \big(\psi^{m,\de},(\de,1-\de)\big),
\\
&\ \ \ \ \ \  \ \ \ E_\infty  \big(\psi^{m,\de},(1-\de,1)\big)\Big\} 
\end{split}
\]
and since $\psi^{m,\de}=\psi^\infty$ on $(\de,1-\de)$, we have
\beq  \label{2.6}
\begin{split}
E_\infty \big(\psi^{m,\de},(0,1)\big)\, 
 \leq & \, \max\bigg\{ E_\infty  \big(\psi^{m,\de},
(0,\de)\big), \,  E_\infty  \big(\psi^\infty,(0,1 )\big),
\\
& \ \ \ \ \ \  \ \ \ E_\infty  \big(\psi^{m,\de},(1-\de,1)\big)\bigg\}.
\end{split}
\eeq
Combining \eqref{2.4}-\eqref{2.6} and \eqref{2.3}, we get
\beq \label{2.7}
\begin{split}
E_\infty \big(u^\infty,(0,1)\big)\, & \leq\, \liminf_{m\ri\infty}\Big( \max\Big\{  E_\infty  \big(\psi^{m,\de},
(0,\de)\big), \,  E_\infty  \big(\psi^\infty,(0,1 )\big),
\\
& \hspace{75pt} E_\infty  \big(\psi^{m,\de},(1-\de,1)\big)\Big\} \Big)
\\
& \leq\, \max\Big\{  E_\infty  \big(\psi^\infty,(0,1 )\big),\,  E_\infty  \big( \psi^{\infty, \de}  ,
(0,\de)\big) ,
\\
&  \hspace{40pt}  E_\infty  \big(  \psi^{\infty,\de},(1-\de,1)\big)\Big\}.
\end{split}
\eeq
Let us now denote the difference quotient of a function $v : \R \larrow \R^N$ as $\D^{1,t}v(x):=\frac{1}{t}[v(x+t)-v(x)]$. Then, we may write 
\[
\ \ \ \ 
\begin{split}
&\text{$\D\psi^{\infty,\de}(x) \,= \,\D^{1,\de}\psi^\infty(0)$, \ \ \ \ \ $x\in (0,\de)$}, 
\\
&\text{$\D\psi^{\infty,\de}(x)\,= \,\D^{1,-\de}\psi^\infty(1)$, \ \ \ $x\in (1-\de,1)$,
}
\end{split}
\]
Note now that
\beq \label{2.8}
\left\{\ \ \
\begin{split}
E_\infty  \big(\psi^{\infty,\de},
(0,\de)\big) \, &=\, \max_{0\leq x\leq \de} \msL\Big(x, \psi^{\infty,\de}(x),\D^{1,\de}\psi^\infty(0)\Big),
\\
E_\infty  \big(\psi^{\infty,\de},
(1-\de,1)\big) \, &=\, \max_{1-\de \leq x\leq 1} \msL\Big(x, \psi^{\infty,\de}(x),\D^{1,-\de}\psi^\infty(1)\Big).
\end{split}
\right.
\eeq
In view of  \eqref{2.7}-\eqref{2.8}, it is suffices to prove that there exist an infinitesimal sequence $(\de_i)_{i=1}^\infty$ such that
\beq \label{2.9}
\begin{split}
E_\infty  \big(\psi^\infty ,
(0,1)\big) \, \geq \, \max\bigg\{ \limsup_{i\ri\infty}  &
\max_{[0,\de_i]}\ \msL\Big(\cdot, \psi^{\infty,\de_i},\D^{1,\de_i}\psi^\infty(0)\Big),
\\
\limsup_{i\ri\infty}  & \max_{[1-\de_i,1]}\ \msL\Big(\cdot, \psi^{\infty,\de_i},\D^{1,-\de_i}\psi^\infty(1)\Big) \bigg\}.
\end{split}
\eeq
The rest of the proof is devoted to establishing \eqref{2.9}. Let us begin by recording for later use that
\beq \label{2.10}
\left\{\ \ \ 
\begin{split}
\max_{0\leq x\leq \de} & \Big| \psi^{\infty,\de}(x) -\psi^\infty(0)\Big| \, \larrow 0, \ \ \text{ as }\de \ri 0,
\\
 \max_{1-\de\leq x\leq 1} & \Big| \psi^{\infty,\de}(x) -\psi^\infty(1)\Big| \, \larrow 0, \ \ \text{ as }\de \ri 0.
\end{split}
\right.
\eeq
Fix a generic $u\in W^{1,\infty}(\Om,\R^N)$, $x\in [0,1]$ and $0<\e<1/3$ and define 
\[
A_\e(x)\,:=\, [x-\e,x+\e] \cap [0,1].
\]
We claim that there exist an increasing modulus of continuity $\om \in C(0,\infty)$ with $\om(0^+)=0$ such that
\beq \label{2.11}
E_\infty\big(u,A_\e(x) \big)\, \geq\, \underset{y \in A_\e(x)}{\ess\,\sup}\ \msL\Big(x,u(x),\D u(y) \Big)\, -\, \om(\e).
\eeq 
Indeed for a.e.\ $y \in A_\e(x)$ we have  $|x-y|\leq \e$ and by the continuity of $\msL$ and the essential boundedness of the derivative $Du$, there exist $\om$ such that 
\[
\Big| \msL\Big(x,u(x),\D u(y) \Big) - \msL\Big(y,u(y),\D u(y) \Big) \Big| \, \leq\, \om(\e)
\]
for a.e.\ $y \in A_\e(x)$, leading directly to \eqref{2.11}. Now, we show that
\beq  \label{2.12}
\sup_{A_\e(x)} \left\{\limsup_{t\ri0}\ \msL\Big(x,u(x),\D^{1,t}u(y) \Big) \right\}\, \leq\,
\underset{A_\e(x)}{\ess\,\sup}\ \msL\Big(x,u(x),\D u(y) \Big).
\eeq
Indeed, for any  Lipschitz function $u$, we have
\beq \label{2.13}
\D^{1,t}u(y)\, =\, \frac{u(y+t)-u(y)}{t} \, =\, \int_0^1 Du(y+\lambda t)\,d\lambda,
\eeq
when $y,y+t\in A_\e(x),\ t\neq0$. Further, for any $x\in \Om$ the function $ \msL(x,u(x),\cdot)$ is level-convex and the Lebesgue measure on $[0,1]$ is a probability measure, thus Jensen's inequality for level-convex functions (see e.g.\ \cite{BJW1, BJW2}) yields
\[
\begin{split}
\msL\Big(x,u(x),\D^{1,t}u(y) \Big)\,
& =\,\msL\left(x,u(x), \int_0^1 Du(y+\lambda t)\,d\lambda\right) 
\\ 
&\leq\,\underset{0\leq\lambda\leq 1}{\ess\,\sup}\ \msL\Big(x,u(x),\D u(y+\lambda t) \Big),
\end{split}
\]
when $y\in A_\e(x)$, $0<x<1$. Consequently,
\[
\begin{split}
\underset{A_\e(x)}{ \sup}\,  \bigg\{\limsup_{t\ri0}& \ \msL\Big(x,u(x),\D^{1,t}u(y) \Big)\bigg\}
\\
&\leq\,\underset{A_\e(x)}{\sup}
\left\{ 
\limsup_{t\ri0}\left[
\underset{0\leq\lambda\leq 1}{\ess\,\sup}\ \msL\Big(x,u(x),\D u(y+\lambda t) 
\right]  \right\}
\\ 
&\leq\,\underset{A_\e(x)}{\sup} \left\{\lim_{s\ri0}
\left[
\underset{y-s\leq z \leq y+s}{\ess\,\sup}\ \msL\Big(x,u(x),\D u(z) \Big)\right]\right\}
\\ 
&=\,\underset{A_\e(x)}{\ess\,\sup}\ \msL\Big(x,u(x),\D u(y) \Big),
\end{split}
\]
as desired. Above we have used the following known property of the $L^\infty$ norm (see e.g.\ \cite{C})
\[
\|f\|_{L^\infty(\Om)}\, =\, \sup_{x\in \Om}\bigg(\lim_{\e\ri 0}\bigg\{\underset{(x-\e,x+\e)}{\ess\,\sup}|f|\bigg\}\bigg).
\]
Note now that by \eqref{2.11} we have
\[
\begin{split}
E_\infty\big(u,(0,1)\big)\, \geq\, E_\infty\big(u,A_\e(x)\big) \, \geq\, \underset{A_\e(x)}{\ess\,\sup}\ \mathscr{L}\left(x,u(x),\D u(y)\right) \, -\, \om(\e)
\\
\end{split}
\]
which combined with \eqref{2.12} leads to
\[
\begin{split}
E_\infty\big(u,(0,1)\big)\, &\geq \, \sup_{A_\e(x)} \left(\limsup_{t\ri0}\ \msL\Big(x,u(x),\D^{1,t}u(y) \Big) \right)\, -\, \om(\e)
\\
&\geq \, \limsup_{t\ri0}\ \left(\msL\Big(x,u(x),\D^{1,t}u(x) \Big) \right)\, -\, \om(\e).
\end{split}
\]
By passing to the limit as $\e\ri 0$ we get
\beq \label{2.14}
\ \ E_\infty\big(u,(0,1)\big)\, \geq \  \limsup_{t\ri0}\ \left(\msL\Big(x,u(x),\D^{1,t}u(x) \Big) \right),
\eeq
for any fixed $u\in W^{1,\infty}(\Om,\R^N)$ and $x\in [0,1]$. Now, since  
\[
\big|\D^{1,t}u(x)\big| \, \leq \, \|\D u\|_{L^\infty(\Om)}, \ \ \ \ x\in (0,1),\ t\neq 0,
\]
for any \emph{finite} set of points $x\in (0,1)$, there is a common infinitesimal sequence $(t_i(x))_{i=1}^\infty$ such that 
\beq \label{2.15}
\text{the limit vectors } \lim_{i\ri \infty}\D^{1,t_i(x)}u(x)\  \text{ exists in }\R^N.
\eeq
Utilising the continuity of $\mathscr{L}$ together with \eqref{2.14}-\eqref{2.15} we obtain
\beq \label{2.16}
\begin{split}
E_\infty\big(u,(0,1)\big)\, &  \geq \  \limsup_{i\ri \infty}\ \msL\Big(x,u(x),\D^{1,t_i(x)}u(x) \Big)
\\
&= \,  \mathscr{L}\left(x,u(x),\lim_{i\ri \infty}\D^{1,t_i(x)}u(x)\right).
\end{split}
\eeq
Now we apply \eqref{2.16} to $u\,=\, \psi^\infty$ and $x\in\{ 0, 1\}$ to deduce the existence of a sequence $(\de_i)_{i=1}^\infty$ along which 
\beq \label{2.17}
\text{ the limit vectors }\ \lim_{i\ri \infty}\D^{1,\de_i}\psi^\infty(0), \ \lim_{i\ri \infty}\D^{1,-\de_i}\psi^\infty(1) \ \ \text{ exist in $\R^N$}
\eeq
and also
\beq \label{2.18}
\begin{split}
E_\infty\big(\psi^\infty,(0,1)\big)\, \geq \, \max\bigg\{ & \mathscr{L}\left(0,\psi^\infty(0),\lim_{i\ri \infty}\D^{1,\de_i}\psi^\infty(0)\right)
,\\
& \mathscr{L}\left(1,\psi^\infty(1),\lim_{i\ri \infty}\D^{1,-\de_i}\psi^\infty(1)\right)
\bigg\}.
\end{split}
\eeq
By recalling \eqref{2.8},  \eqref{2.10} and \eqref{2.17}, for $\de=\de_i$ we obtain  
\beq \label{2.19}
\begin{split}
\lim_{i\ri\infty} E_\infty  \big(\psi^{\infty,\de_i},(0,\de_i)\big) \,
&=\, \lim_{i\ri\infty} \max_{[0,\de_i]}\, \msL\Big(\cdot, \psi^{\infty,\de_i},\D^{1,\de_i}\psi^\infty(0)\Big)
\\
&=\, \mathscr{L}\left(0,\psi^\infty(0),\lim_{i\ri \infty}\D^{1,\de_i}\psi^\infty(0)\right),
\end{split}
\eeq
and also
\beq \label{2.20}
\begin{split}
\lim_{i\ri\infty} E_\infty  \big(\psi^{\infty,\de_i},(1-\de_i,1)\big) \, 
&=\, \lim_{i\ri\infty} \max_{[1-\de_i,1]} \msL\Big(\cdot, \psi^{\infty,\de_i},\D^{1,-\de_i}\psi^\infty(1)\Big)
\\
&=\, \mathscr{L}\left(1,\psi^\infty(1),\lim_{i\ri \infty}\D^{1,-\de_i}\psi^\infty(1)\right).
\end{split}
\eeq
By putting together \eqref{2.18}-\eqref{2.20}, \eqref{2.9} ensues and we conclude the proof.    \qed \ms

\ms

\ms

\end{document}